\documentclass{amsart}
\usepackage{color}

\usepackage{amsmath} 
\usepackage{amssymb}
\usepackage{mathrsfs}
\usepackage{mathtools}

\newtheorem{theorem}{Theorem}[section] 
\newtheorem{claim}[theorem]{Claim}

\newtheorem{conclusion}[theorem]{Conclusion}
\newtheorem{observation}[theorem]{Observation}

\theoremstyle{definition}
\newtheorem{definition}[theorem]{Definition}
\newtheorem{dc}[theorem]{Definition/Claim}

\newtheorem{fact}[theorem]{Fact}

\newtheorem{discussion}[theorem]{Discussion}

\theoremstyle{remark}
\newtheorem{remark}[theorem]{Remark}
\newtheorem{concr}[theorem]{Concluding Remarks}
\newtheorem{question}[theorem]{Question}
\newtheorem{notation}[theorem]{Notation}

\newcommand{\Th}{{\rm Th}}

\newcommand{\otp}{{\rm otp}}

\newcommand{\lfgr}{{\rm lfgr}}

\newcommand{\olive}{{\rm olive}}
\newcommand{\Univ}{{\rm Univ}}

\newcommand{\Sym}{{\rm Sym}}

\newcommand{\qf}{{\rm qf}}
\newcommand{\Qr}{{\rm Qr}}

\newcommand{\JEP}{{\rm JEP}}
\newcommand{\INV}{{\rm INV}}

\newcommand{\nacc}{{\rm nacc}}

\newcommand{\bd}{{\rm bd}}

\newcommand{\gr}{{\rm gr}}

\newcommand{\univ}{{\rm univ}}

\newcommand{\LST}{{\rm LST}}
\newcommand{\Mod}{{\rm Mod}}

\newcommand{\NSOP}{{\rm NSOP}}
\newcommand{\SOP}{{\rm SOP}}

\newcommand{\grp}{{\rm grp}}

\newcommand{\Rang}{{\rm Rang}}

\newcommand{\rest}{{\restriction}}
\newcommand{\dom}{{\rm dom}}

\newcommand{\wilog}{{\rm without loss of generality}}
\newcommand{\Wilog}{{\rm Without loss of generality}}

\newcommand{\then}{{\underline{then}}}
\newcommand{\when}{{\underline{when}}}
\newcommand{\Then}{{\underline{Then}}}

\newcommand{\If}{{\underline{if}}}
\newcommand{\Iff}{{\underline{iff}}}

\newcommand{\cA}{{\mathscr A}}

\newcommand{\gA}{{\mathfrak A}}
\newcommand{\gC}{{\mathfrak C}}

\newcommand{\gd}{{\mathfrak d\/}}

\newcommand{\cF}{{\mathscr F}}

\newcommand{\bbL}{{\mathbb L}}

\newcommand{\cM}{{\mathscr M}}

\newcommand{\gk}{{\mathfrak k}}

\newcommand{\bbP}{{\mathbb P}}
\newcommand{\cP}{{\mathscr P}}

\newcommand{\bbQ}{{\mathbb Q}}
\newcommand{\varp}{{\varepsilon}}
\newcommand{\bbZ}{{\mathbb Z}}

\newcommand{\cW}{{\mathscr W}}

\newcommand{\cf}{{\rm cf}}

\newcount\skewfactor
\def\mathunderaccent#1#2 {\let\theaccent#1\skewfactor#2
\mathpalette\putaccentunder}
\def\putaccentunder#1#2{\oalign{$#1#2$\crcr\hidewidth
\vbox to.2ex{\hbox{$#1\skew\skewfactor\theaccent{}$}\vss}\hidewidth}}

\newenvironment{PROOF}[2][\proofname.]
   {\begin{proof}[#1]\renewcommand{\qedsymbol}{\textsquare$_{\rm #2}$}}
   {\end{proof}}

\usepackage[hidelinks]{hyperref}

\begin{document}
\makeatletter\def\shfiuwefootnote{\gdef\@thefnmark{}\@footnotetext}\makeatother\shfiuwefootnote{Version 2024-01-24\_2. See \url{https://shelah.logic.at/papers/1029/} for possible updates.}

\title {No universal group in a cardinal \\ Sh1029}

\author{Saharon Shelah}

\address{Einstein Institute of Mathematics\\
Edmond J. Safra Campus, Givat Ram\\
The Hebrew University of Jerusalem\\
Jerusalem, 91904, Israel\\
and \\
Department of Mathematics\\
Hill Center - Busch Campus \\ 
Rutgers, The State University of New Jersey \\
110 Frelinghuysen Road \\
Piscataway, NJ 08854-8019 USA}
 
\email{shelah@math.huji.ac.il}

\urladdr{http://shelah.logic.at}

% \thanks{The author thanks the Israel Science Foundation for partial
% support of this work.  Grant No. 1053/11, Publication 1029.  
% In early versions (up to 2017) the author thanks Alice Leonhardt for the beautiful 
% typing. First typed: March 27, 2013.}

\thanks{First typed: March 27, 2013. The author would like to thank the Israel Science Foundation for partial
support of this work.  Grant No. 1053/11. In early versions (up to 2017) the author thanks Alice Leonhardt for the beautiful  typing. In later versions, the author would like to thank the typist for his work and is also grateful for the generous funding of typing services donated by a person who wishes to remain anonymous.  References like [Sh:950, Th0.2=Ly5] mean that the internal label of Th0.2 is y5 in Sh:950. The reader should note that the version in the author's website is usually more up-to-date than the one in arXiv. This is publication number 1029 in Saharon Shelah's list.}



% Previous version: November 17, 2015 

%(2013/June/21; was f1320; copied and revised to sh1029 on July 7, 2013)

\date{May 15, 2023}

\begin{abstract}
    For many classes of models, there are universal members in any cardinal
    $\lambda$ which ``essentially satisfies GCH, i.e. $\lambda =
    2^{< \lambda}$", in particular for the class of a complete first order
    $T$ (well, if at least if $\lambda > |T|$).   But if the class is ``complicated enough", e.g. the
    class of linear orders, we know that if $\lambda$ is ``regular and not
    so close to satisfying GCH" then there is no universal member.  Here we find new sufficient conditions (which we call the olive property), not covered by earlier cases (i.e. fail the so-called $\rm SOP_4$).  The
    advantage of those conditions is witnessed by proving that the class
    of groups satisfies one of those conditions.

    This version has minor changes , compared to the earlier one.
\end{abstract}

\maketitle
\numberwithin{equation}{section}
\setcounter{section}{-1}

%%%%%%%%%%%%%%%%%%%%%%%%%%%%%%%%%%%%%
% Annotated content
\newpage

\centerline{Annotated Content}

\S0 \quad Introduction, (labels y,z), pg.\pageref{0}

\S1 \quad The olive property, (label d), pg.\pageref{1}

\begin{enumerate}
    \item[${{}}$]  [We give definitions of some versions of the olive property and give an example failing the $\SOP_4$.  We phrase relevant set theoretic conditions like $\Qr_1$ (slightly weaker than those used earlier).  Then we give complete proof using $\Qr_1(\chi_2,\chi_1,\lambda)$ to deduce $\Univ(\chi_1,\lambda,\gk)  \ge \chi_2$ so there is no universal in the class $\gk$ in the cardinal $\lambda$, when $\gk$ has the olive property.]
\end{enumerate}

\S2 \quad The class of groups have the olive property, (label s), 
pg.\pageref{2}

\begin{enumerate}
    \item[${{}}$]  [We prove the stated result.  We also deal with the non-existence of universal structures for pairs of classes, e.g.  the pair (locally finite groups, groups).]
\end{enumerate}

%%%%%%%%%%%%%%%%%%%%%%%%%%%%%%%%%%%%%%%
% Introduction 
\newpage

\section{Introduction}\label{0}

\subsection{Background and open questions}\label{0A}\

A natural and old question is how to characterize the class of cardinals $\lambda$ in which a class $\mathbf k$ has a universal member, where $\gk$ is, e.g. the class of models of a first-order theory $T$ with elementary embeddings.

On history see Kojman-Shelah \cite{Sh:409} and later Dzamonja \cite{Dj05}.  Recall that if $\lambda = 2^{< \lambda} > \aleph_0$ then many classes have a universal member in $\lambda$, so assuming GCH, we know when there is a universal model in every $\lambda > |T|$.

For transparency, we consider a first-order countable $T$. Now the class $\gk_T = (\Mod_T,\prec)$ of models of $T$ with elementary
embeddings have the amalgamation property, the JEP, and satisfies: $A \subseteq M \in \Mod_T \Rightarrow (\exists N \in \Mod_T)(A \subseteq N \prec M
\wedge \|N\| = |A| + \aleph_0)$.  From such classes (or just a.e.c. with amalgamation and JEP with $\LST(\gk)$ instead of $\aleph_0$) it follows that if $\lambda = \lambda^{< \lambda} > \aleph_0$, \then \, there is a saturated (or universal homogeneous) 
$M \in \Mod_T$ of cardinality $\lambda$ which implies it is universal for $\gk_T$.  If $\lambda = 2^{<\lambda} > \aleph_0$ does not satisfy $\lambda = \lambda^{< \lambda}$, still there is a so-called special model which is universal.  So we are interested in the cases where G.C.H. fail.

Recall that on the one hand, Kojman-Shelah \cite{Sh:409} shows that if $T$ is the theory of dense linear orders or just $T$ has the strict
order property, then $T$ fails (in a strong way) to have a universal member in regular cardinals in which cardinal arithmetic is ``not close to GCH";
(for regular $\lambda$ this means there is a regular $\mu$ such that $\mu^+ < \lambda < 2^\mu$, while for singular $\lambda$ we need of course $\lambda < 2^{< \lambda}$ and a very weak pcf condition).

By \cite{Sh:500}, we can weaken ``the strict order property" to the 4-strong order property $\SOP_4$.  On consistency see Dzamonja-Shelah \cite{Sh:614}. 

Natural questions (we shall address some of them) are the following:

\begin{question}\label{y2}
    1) Is there a weaker condition (on $T$) than $\SOP_4$ which suffices?
    
    2) Can we find a best one?
    
    3) Can we find such a condition satisfied for some theory $T$ which is  $\NSOP_3$?
\end{question}

\begin{question}\label{y5}
    1) Is there $T$ with the class $\Univ(T) \setminus (2^{\aleph_0})^+$ strictly
    smaller than the one for linear orderings, see Definition \ref{z8}(2);  is it better if we restrict ourselves to regular cardinals above $2^{\aleph_0}$?  
    
    2) Can we get the above to be $\{\lambda \colon \lambda = 2^{< \lambda} \}$?
    
    3) What about singular cardinals?
\end{question}

\begin{question}\label{y7}
    1) Is it consistent that the class of linear orders has a universal member in $\lambda$ such that $2^{< \lambda} > \lambda > 
    2^{\aleph_0}$.  (For $\lambda = \aleph_1 < 2^{\aleph_0}$, the answer is yes, see \cite{Sh:100}, a  more detailed version is in preparation). 
    
    2) The same can be asked for some theory with $\SOP_4$ or the olive property (defined below, e.g. in Definition \ref{y36}.
\end{question}

Recall that by Shelah-Usvyatsov \cite{Sh:789} the class of groups  has $\NSOP_4$ but has $\SOP_3$, so it was not clear where it stands.

\begin{question}\label{y8}
    1) Where does the class of groups stand (concerning the existence of a universal member in a cardinal)?
    
    2) Is it consistent that there is a universal locally finite group of cardinality $\aleph_1$? of cardinality $\beth_\omega$? $\lambda = \beth^+_\omega$? of regular cardinals $\lambda \in (\beth^+_\omega,\beth_{\omega +1})$? of other cardinals $\lambda  < \lambda^{\aleph_0}$?
    
    Recall (Grossberg-Shelah \cite{Sh:174}) that if 
    $\mu$ is strong limit of cofinality $\aleph_0$ 
    above a a compact cardinal, \then \, there is a universal locally finite group of cardinality $\mu$ but if $\mu = \mu^{\aleph_0}$ then there is no one.
\end{question}

Concerning singulars

\begin{question}\label{y11}
    Does $\theta = \cf(\theta)$ and $\theta^{+2} < \cf(\lambda) < \lambda < 2^\theta$ implies $\lambda < \univ(\lambda,T)$?
\end{question}

\begin{question}\label{y14}
    0) Characterize the failure of the criterion of \cite{Sh:457}, D\v{z}amonja-Shelah \cite{Sh:614}(for consistency).
    
    1) Does $\SOP_3$ (or something weaker) suffice for no universal in $\lambda$ when $\mu = \mu^{<\mu} \ll \lambda < 2^\mu$?
    
    2) Which theories $T$ fail to have a universal 
    in $\lambda$ when $\lambda = \mu^{++} = 2^\mu < 2^{\mu^+}$?
    
    3) We may consider weaker properties 
    of $T$ for no universal in $\lambda,\mu = \mu^{< \mu} \ll \lambda < 2^\mu$. 
    
    4) Sort out the variants of the olive property (defined below, e.g. in Definition \ref{y36}). 
\end{question}

\begin{discussion}\label{y20}
    Even for linear orders, the case
    
    \begin{enumerate}
    \item[$(*)^1_\lambda$]  \underline{successor case}:  $\lambda = \mu^+,\lambda < 2^\mu$ and $2^{< \mu} \le \lambda$ (e.g.
     for transparency $\mu = \mu^{< \mu}$) is not 
    resolved as we do not necessarily have $\bar C = \langle C_\delta \colon \delta \in S^\lambda_\mu\rangle$ guessing clubs,
    recall that by \cite{Sh:409} if $2^\theta > \lambda = \cf(\lambda) >
    \theta = \cf(\theta),2^\theta \le \univ(\lambda,T)$, so if $\lambda$
    is a successor cardinal $> \lambda_n$, the only open case is $(*)^2_\lambda$.  Similarly, if $\lambda$ is a limit cardinal the only open case is
    
    \item[$(*)^2_\lambda$]  \underline{limit case}: $\lambda$ is singular.
    \end{enumerate}
    
    In this case, there are strong pcf restrictions (see \cite{Sh:409}), so advancement there may eliminate the case.
    
    By some earlier results (see \cite{Sh:457})  if $\mu = 2^\kappa$, (so $\mu$ is not a strong limit cardinal), and there is no universal in $\lambda$ \then \, there was 
    a sequence $\langle \Lambda_\delta \colon \delta \in
    S^\lambda_\mu\rangle,\Lambda_\delta \subseteq {}^{(C_\delta)}\mu$ of
    cardinality $\lambda$ such that for every sequence $\langle \eta_\delta
    \subseteq {}^{(C_\delta)}\mu \colon \delta \in S^\lambda_\mu\rangle$ there is
    club $E$ of $\lambda$ such that for every $\delta \in E \cap
    S^\lambda_\mu$ there is a 
    $\nu \in \Lambda_\delta$ such that the functions $\eta_\delta,\nu$ agree on $E \cap \nacc(C_\delta)$.
    
    Using a more complicated $T$ we can replace ${}^{C_\delta}\mu$ by ${}^{(C_\delta \times D_\delta)}\mu$ so the agreement above is on the
    product $(E \cap \nacc(C_\delta)) \times (E  \cap \nacc(C_\delta))$  but of unclear value.  
    
    On subsequent works and more on consistency, see \cite{Sh:F1330} and \cite{Sh:F1414}.
    
    We thank the referee and Thilo Weinert for their helpful comments.
\end{discussion}

\subsection{What is accomplished}\

What do we achieve?  We introduce the ``olive property" which is a sufficient condition for a class to have a universal  member in $\lambda$ only if $\lambda$ is ``close to satisfying G.C.H.", similar to the linear order case.   This condition is weaker than $\SOP_4$, hence giving a positive answer to Question \ref{y2}(1).  But the
condition implies $\SOP_3$ so it does not answer Question \ref{y2}(3), also it
is totally unclear whether it is best in any sense and whether its negation has interesting consequences.

However, it answers Question \ref{y8}(1) to a large extent because the class of
groups have the olive property and we can also deal with locally
finite groups answering some cases of Question \ref{y8}(2); see \S2.  
Also, we try to formalize
conditions sufficient for non-existence, see Definition \ref{d13}.  
As the reader may find the definition of the (variants of the) olive property
opaque, we define a simple case used for the class of groups, and  the reader may then look first at the class of groups in \S2.

\begin{definition}\label{y36}
    A (first order) universal theory $T$ has the \emph{olive property} \when \, there are $(\varphi_0,\varphi_1,\psi)$ and a model $\gC$ of $T$ such that:
    
    \begin{enumerate}
        \item[(a)]  for some $m,\varphi_0 = \varphi_0(\bar x_{[m]}, \bar
        y_{[m]}),\varphi_1 = \varphi_1(\bar x_{[m]},\bar y_{[m]}),\psi =
        \psi(\bar x_{[m]},\bar y_{[m]},\bar z_{[m]})$ are quantifier free
        formulas (and $\bar x_{[m]},\bar y_{[m]},\bar z_{[m]}$ are $m$-tuples
        of variables, see Notation \ref{z3} below)
        
        \item[(b)]  for every $k$ and $\bar f = \langle f_\alpha \colon \alpha < k 
        \rangle$ where $f_\alpha$ is a function from $\alpha$ to $\{0,1\}$ we can find a model $M$ of $T$ and  $\bar a_\alpha \in {}^m M$ for $\alpha < k$ such that:
        
        \begin{enumerate}
            \item[$(\alpha)$]  $\varphi_\iota[\bar a_\alpha,\bar a_\beta]$ for
            $\alpha < \beta < k$ when $\iota = f_\beta(\alpha),$
            
            \item[$(\beta)$]  $\psi[\bar a_\alpha,\bar a_\beta,\bar a_\gamma]$ when $\alpha < \beta < \lambda$ and $f_\gamma \rest [\alpha,\beta]$ is
            constantly 0.
        \end{enumerate}
        
        \item[(c)] there are no $\bar a_\ell \in {}^m M$ for
        $\ell=0,1,2,3$ such that the following conditions are\footnote{in the class of
        groups, in clause $(\alpha),\varphi_0[\bar a_0,\bar a_1], \varphi_1[a_1,a_2],\varphi_1[a_1,\bar a_3]$ suffice.} satisfied in $M$:
        
        \begin{enumerate}
            \item[$(\alpha)$]  $\varphi_0[\bar a_0,\bar a_\ell]$ for
            $\ell=1,2,3,\varphi_1[\bar a_1,\bar a_\ell]$ for $\ell=1,2$ and $\varphi_0[\bar a_2,\bar a_3],$
            
            \item[$(\beta)$]  $\psi[\bar a_0, \bar a_2, \bar a_3].$
        \end{enumerate}
    \end{enumerate}
\end{definition}

\begin{concr}
    Concerning some things not addressed here we note the following.
    
    1) Concerning the proof here of ``there is no universal" we can carry it via defining invariants parallel to Kojman-Shelah \cite{Sh:409} such that (for transparency $\lambda$ is regular uncountable, see
    Definition \ref{z6}(5), (7)):
     
    \begin{enumerate}
        \item[$(\ast)$]
    
        \begin{enumerate}
            \item[(a)] if $M \in \Mod_{T,\lambda}$ then $\INV_\lambda(M)$ is a set of cardinality $\le \lambda$ or just  $\le \chi < 2^\lambda,$
    
            \item[(b)] if $M_1,M_2 \in \Mod_{T,\lambda}$ and $M_1$ is elementarily embeddable into $M_2$ then  $\INV_\lambda(M_1) \subseteq \INV_\lambda(M_2),$
    
            \item[(c)] there is a set of $2^\lambda$ objects $\mathbf x$ such that $(\exists M \in \Mod_{T,\lambda})(\mathbf x \in 
            \INV_\lambda(M))$.
        \end{enumerate} 
    \end{enumerate}
    
    2) We can use more complicated versions of the olive property.  In the proof we use one $\delta$ and then one $\alpha \in \nacc(C_\delta)
    \cap E$ (or less), but we may use several ordinals $\alpha$ resulting in more complicated versions.  This will become more pressing if we have a complimentary property, guaranteeing ``no universal" or some variant.
\end{concr}

\subsection{Preliminaries}\label{0C}\

\begin{notation}\label{z3}
    1) Let $\bar x_{[I]} = \langle x_t \colon t \in I \rangle$ and similarly $\bar y_{[I]},\bar x_{[I],\alpha}$, etc. where $\bar x_{[I],\ell} =
    \langle x_{t,\ell} \colon t \in I \rangle$.
    
    2) For a first-order complete $T,\gC_T$ is the ``monster model of $T$" omitting $T$ if clear from the context.
\end{notation}

\begin{definition}\label{z6}
    1) For a set $A,|A|$ is its cardinality but for a structure $M$ its cardinality is $\|M\|$ while its universe is $|M|$; this applies e.g. to groups.
    
    2) We use $G,H$ for groups, $M,N$ for general models.
    
    3) Let $\gk$ denote a pair $(K_{\gk},\le_{\gk})$, or we may say a class (of
    models) $\gk$, where:
    
    \begin{enumerate}
        \item[(a)]  $K_{\gk}$ is a class of $\tau_{\gk}$-structures where $\tau_{\gk}$ is a vocabulary,
        
        \item[(b)]  $\le_{\gk}$ is a partial order on $K_{\gk}$ such that $M \le_{\gk} N \Rightarrow M \subseteq N,$
        
        \item[(c)]  both $K_{\gk}$ and $\le_{\gk}$ are closed under isomorphisms.
    \end{enumerate}
    
    \begin{enumerate}
        \item[4) \ \ (a)]  We may write only $K_{\gk}$
          \when \, $\le_{\gk}$ is being a submodel,
        
        \item[(b)]  We say $f \colon M \rightarrow N$ is a  $\le_{\gk}$-embedding \when \, $f$ is
        an isomorphism from $M$ onto some $M_1 \le_{\gk} N$.
    \end{enumerate}
    
    5) If $T$ is a first-order theory then $\Mod_T$ is the pair $(\!\!\!\mod_T,\le_T)$ where $\text{ mod}_T$ is the class of models of $T$ and $\le_T$ is: $\prec$ if $T$ is complete, $\subseteq$ if $T$ is not complete.
    
    6) We may write $T$ instead of $\Mod_T$, e.g. in Definition \ref{z8} below.
    
    7) For a class $K$ of structures $K_\lambda = \{M \in K \colon \|M\| = \lambda\}$.
\end{definition}

\begin{definition}\label{z8}
    1) For a class $\gk$ and a cardinal $\lambda$, a set $\{M_i \colon i<i^*\}$ of models from $K_{\gk}$, is
    \emph{jointly $(\lambda,\gk)$-universal} \when \, for 
    every $N \in K_{\gk}$ of size $\lambda$,
    there is an $i<i^*$ and an $\le_{\gk}$-embedding of $N$ into $M_i$.
    
    2) For $\gk$ and $\lambda$ as above, let (if $\mu = \lambda$ we may
    omit $\mu$)
    
    \begin{equation*}
        \begin{array}{clcr}
        \univ(\lambda,\mu,\gk) \coloneqq &\min\{|\cM| \colon \cM \text{ is a family of
          members of } K_{\gk} \text{ each}\\
          &\text{ of cardinality } \le \mu \text{ which is jointly } \gk\text{-universal for } \lambda\}
        \end{array}
    \end{equation*}
    
    Let $\Univ(\gk)  \coloneqq \{\lambda \colon \univ(\lambda,\gk)=1\}$.
    
    3) For a pair\footnote{This will be used for $\gk_1$ being the class of locally finite groups and $\gk_2$ being the class of groups.}
    $\bar{\gk} = (\gk_1,\gk_2)$ of classes with $\gk_\iota = (K_{\gk_\iota},\le_{\gk_\iota})$ as in Definition \ref{z6}(3) for $\iota=1,2$ such that $\tau(\gk_1) = \tau(\gk_2)$ and  $K_{\gk_1} \subseteq K_{\gk_2}$, let $\univ(\lambda,\mu,\bar{\gk})$ be the minimal $|\cM|$ such
    that $\cM$ is a family of members of $K_{\gk_2}$ each of cardinality
    $\le \mu$ such that every $M \in K_{\gk_1}$ of cardinality $\lambda$ can be $\le_{\gk_2}$-embedded into some member of $\cM$.
\end{definition}

Dealing with a.e.c.'s (see \cite{Sh:h}) we have the following:

\begin{definition}\label{z10}
    1) We say that a formula $\varphi = \varphi(\bar x_{[I]})$, in any logic, is $\gk$-upward preserved \when \, $\tau_\varphi \subseteq \tau_{\gk}$ and if $M \le_{\gk} N$ and $\bar a \in {}^I M$ 
    then $M \models \varphi[\bar a]$ implies $N \models \varphi[\bar a]$.
    
    2) For $\bar{\gk}$ as in Definition \ref{z8}(3) we say that  a pair $\bar\varphi(\bar x_{[I]}) = 
    (\varphi_1(\bar x_{[I]}),  \varphi_2(\bar x_{[I]}))$ is
    $\bar{\gk}$-upward preserving \when \, $\tau_{\varphi_1} \cup
    \tau_{\varphi_2} \subseteq \tau_{\gk_\iota}$ and if $M_\iota \in
    K_{\gk_\iota}$ for $\iota =1,2,\bar a \in {}^I(M_1)$ and $M_1
    \le_{\gk_2} M_2$ then $M_1 \models \varphi_1[\bar a]$ implies $M_2
    \models \varphi_1[\bar a]$.
    
    3) In part (2), if $\varphi_0 = \varphi_1$ then we may write $\varphi$ instead of $\bar\varphi$.  Saying\footnote{Pedantically, a pair is a sequence of length 2 so the Definition \ref{z10}(2), \ref{z10}(3) are incompatible, but the intention should be clear from the context.} that a sequence $\bar\psi$ is $\gk$-upward preserving means that every formula appearing in $\bar\psi$ is $\gk$-upward preserving.
\end{definition}

\begin{definition}\label{z20}
    1) For an ideal $J$ on a set $A$ and a set $B$ let $\mathbf U_J(B) = 
    \min\{|{\cP}| \colon {\cP}$ is a family of subsets of $B$, each of 
    cardinality $\le |A|$ such that 
    for every function $f$ from $A$ to $B$ for some  $u \in {\cP}$ we have $\{a \in A \colon f(a) \in u\} \in J^+\}$.
    
    2) For an ideal $J$ on a set $A$, a cardinal $\theta$ and a set $B$ let
    $\mathbf U^\theta_J(B) = \min \{|\cP| \colon \cP \subseteq [B]^{\le |A|}$ and
    if $f \in {}^A({}^\theta B)$ then for some $u \in \cP$ we have $\{a
    \in A \colon \Rang(f(a)) \subseteq u\} \in J^+\}$.  So $\mathbf U^\theta_J(B) \le \mathbf U_J(|N|^\theta)$.
    
    3) Clearly only $|B|$ matters so we normally write $\mathbf U_J(\lambda)$, (see on it \cite{Sh:589}).
\end{definition}

%%%%%%%%%%%%%%%%%%%%%%%%%%%%%%%%%%%%%%%
% Section 1
\newpage

\section{The Olive property}\label{1}
 
\begin{definition}\label{d2}
    1) (Convention)
    
    \begin{enumerate}
        \item[(a)]  Let $T$ be a first-order theory and $\gC \coloneqq \gC_T$ a monster for $T$,
        
        \item[(b)]
    
        \begin{enumerate}
            \item[$(\alpha)$] $ \Delta \subseteq \bbL(\tau_T)$  a set of formulas,
    
            \item[$(\beta)$] omitting $\Delta$ means $\Delta = \bbL(\tau_T)$ if $T$ is complete, $\Delta$ = set of quantifiers free formula otherwise, and we may write $\qf$ instead of $\Delta.$
        \end{enumerate}
        
        \item[(c)]
        
        \begin{enumerate}
            \item[$(\alpha)$] $ m$ and $n \ge k_\iota \ge 2$ for $\iota=0,1, n \ge k_0 + k_1 \ge 3,\eta \in {}^n 2$ be such that   $\eta(0) = 0$ and $\eta^{-1}\{0\}$ is not an initial segment and $\eta^{-1}\{\iota\}$ has $\ge k_\iota$ members for $\iota=0,1,$
    
            \item[$(\beta)$]  If $\eta(\ell) = \ell \mod 2$ for $\ell < k$ we may write $n$ instead of $\eta.$
        \end{enumerate}
    
        \item[(d)] 
    
        \begin{enumerate}
            \item[$(\alpha)$] If $\bar k = 
            (k_0,k_1),k_1 \le k_0 + 1 \le k_1+1$ we may write   $k_0 + k_1$ instead of $\bar k$ and let $k(\iota) = k_\iota$ for $\iota=0,1,$ 
    
            \item[$(\beta)$] omitting $m$ means ``for some $m$",
    
            \item[$(\gamma)$] omitting $n,\eta,\bar k$ means $n=3,\eta =  \langle 0,1,0\rangle,\bar k = (2,1)$) for some $m$.
        \end{enumerate}
        
        \item[(e)]
        
        \begin{enumerate}
            \item[$(\alpha)$] below we may write  $\psi_\iota = \psi_{\iota,k_\iota}$ and $\varphi_\iota = \psi_{\iota,1}$ for $\iota=0,1$,
    
            \item[$(\beta)$] if $\varphi_0 = \varphi_1 = \varphi$ we may write $\varphi$,
    
            \item[$(\gamma)$]  we may omit $\psi_{3,k}$ when it is a 
            logically true formula.
        \end{enumerate}
    \end{enumerate}
    
    2) We say $T$ has the \emph{$(\Delta,\eta,\bar k,m)$-olive property} \when \, there is a pair $(\bar\psi_0,\bar\psi_1)$ of sequences of formulas from $\Delta$ witnessing it, see (3).
    
    3) We say $(\bar\psi_0,\bar\psi_1)$ \emph{witnesses the $(\Delta,\eta,\bar
    k,m)$-olive property} (for $T$, with the convention above) \when \, (it
    is the case for every $\lambda$, but in this definition, by
    compactness,  $\lambda = \aleph_0$ is enough):
    
    \begin{enumerate}
        \item[(a)]  $\bar\psi_\iota = \langle \psi_{\iota,k}(\bar
        x_{[m],0}, \dotsc,\bar x_{[m],k}) \colon  k = 1,\dotsc,k_\iota
        \rangle$ for $\iota = 0,1$ with $\psi_{\iota,k} \in \Delta,$
        
        \item[(b)$_\lambda$]  for every $\bar f = \langle f_\alpha \colon \alpha < \lambda\rangle$ where $f_\alpha$ is a function from $\alpha$ to $\{0,1\}$, we can find a model $M$ of $T$ and  $\bar a_\alpha \in {}^m M$ for $\alpha < \lambda$
        such\footnote{Actually clause $(\alpha)$ is a specific case of
        clause $(\beta)$ provided that in clause $(\beta)$ we allow $k=1$.
        Similarly for clauses $(c)(\alpha), (\beta)$.} that:
        
        \begin{enumerate}
            \item[$(\alpha)$]  $\varphi_\iota[\bar a_\alpha,\bar a_\beta]$ for
            $\alpha < \beta <\lambda$ when $\iota = f_\beta(\alpha)$, recalling Definition \ref{d2}$(1)(e)(\alpha),$
            
            \item[$(\beta)$]  $\psi_{\iota,k}(\bar a_{\alpha_0},\dotsc,
            \bar a_{\alpha_{k-1}},\bar a_\beta)$ when $k \in
            \{2,\dotsc,k_\iota\}$ and $\alpha_0 < \ldots < \alpha_{k-1} < \beta  <\lambda$ and $f_\beta \rest [\alpha_0,\alpha_{k-1}]$ is constantly $\iota$, so when $k=1$, it holds trivially,
        \end{enumerate}
        
        \item[(c)]  there are no $\bar a_\ell \in {}^m{\gC}$ for $\ell <
        n+1$ such that:
        
        \begin{enumerate}
            \item[$(\alpha)$]  $\varphi_\iota[\bar a_i,\bar a_j]$ for $i<j < n+1$
            and $\eta(i) = \iota,$
            
            \item[$(\beta)$]  if $\iota \in \{0,1\}, $ \, $ k \in \{2,\dotsc,k_\iota\},$ 
            $\ell_0 < \ldots < \ell_{k-1}$ are from $\{\ell < n \colon \eta(\ell)  = \iota\}$ and $\ell_{k-1} < \ell \le n,$ then 
            $\psi_{\iota,k}[\bar a_{\ell_0},\dotsc, \bar a_{\ell_{k-1}}, \bar a_\ell]$.
        \end{enumerate}
    \end{enumerate}
\end{definition}

\begin{remark}\label{d3}
    This fits the classification of properties of such $T$ in \cite[5.15-5.23]{Sh:702}.
\end{remark}

\begin{definition}\label{d4}
    1) Let $K$ be a universal class of $\tau$-models, see Definition \ref{z6}(4).   We say $K$ has  
    the \emph{$\lambda-(\eta,\bar k,m)$-olive property} \when \, some quantifier-free $(\bar\psi_0,\bar\psi_1)$ 
    witnessing it, that is, $(a) + (b)_\lambda  + (c)$ holds (replacing $\gC_T$ by ``in some $M \in K_\lambda$").
    
    2) We say that a class (e.g. an a.e.c.) $\gk = (K_{\gk},\le_{\gk})$ has the \emph{$\lambda-(\eta,\bar k,m)$-property} \when 
    there are $\bar\psi_0,\bar\psi_1$ which are $\gk$-upward preserved formulas in any logic (see Definition \ref{z10})
    and $(a) + (b)_\lambda + (c)$ of Definition \ref{d2} holds, 
    replacing $M$ by ``some $\gC \in K_{\gk}$ of cardinality $\lambda$".
\end{definition}

\begin{remark}\label{d5}
    1) Note that for $T$ first order complete, $\gk = \Mod_T = (\text{mod}_T,\prec)$, Definition \ref{d4}(2) gives Definition \ref{d2} and for $T$ first order universal not complete, $\gk = \Mod_T = (\text{mod}_T,\subseteq)$, Definition \ref{d4}(2) gives Definition \ref{d2}.  Similarly for Definition \ref{d4}(1).
    
    2) Of course, for $T$ first order, the $\lambda$ does not matter.
\end{remark}

\begin{claim}\label{d7}
    Assume $n \ge k_0 + k_1 \ge 3,\eta \in {}^n 2$ and $|\eta^{-1}\{\iota\}| \ge k_\iota \ge 1$ for $\iota = 0,1$. \Then \, there is a complete first-order countable $T$ having the $(\eta,\bar k,1)$-olive property but $T$ is
    $\NSOP_4$, is $\SOP_3$ and is categorical in $\aleph_0$.
\end{claim}

\begin{PROOF}{\ref{d7}}  
    Let $\tau = \{P,Q_0,Q_1\}$ where $P$ is a binary predicate and
    $Q_\iota$ is a $(k_\iota +1)$-place predicate.  Let $T^0_{\eta,\bar k}$ 
    be the following universal theory in $\bbL(\tau)$:
    
    \begin{enumerate}
        \item[$(*)_1$]  a $\tau$-model $M$ is a model of $T^0_{\eta,\bar k}$  
        \Iff \, we cannot embed $N^*_{\eta,\bar k}$ into $M$ where
    
        \begin{enumerate}
            \item[$\oplus$]  $N^*_{\eta,\bar k}$ is the $\tau$-model with universe $\{a_0,\dotsc,a_n\}$ as in $(c)(\alpha), (\beta)$ from   Definition \ref{d2}(3) for $\varphi(x_0,x_1) = P(x_0,x_1), $ and $\psi_\iota(x_0,\dotsc,x_{k(\iota)}) =Q_\iota(x_0,\dotsc,x_{k(\iota)}),$ recalling Definition \ref{d2}$(1)(e)(\gamma)$ and $\ell < k \le n \Rightarrow a_\ell \ne a_k$.
        \end{enumerate}
    \end{enumerate}
    
    Now, 
    
    \begin{enumerate}
        \item[$(*)_2$]  $T^0_{\eta,\bar k}$ has the JEP and the 
        amalgamation property by disjoint union.
    \end{enumerate}
    
    [Why?  Assume that $M_0 \subseteq M_1,M_0 \subseteq M_2$ are models of
    $T_0$ (but abusing notation we allow $M_0$ to be empty) and $|M_1| \cap |M_2| = |M_0|$, we define $M = M_1 \cup M_2,$ that is, 

    \begin{enumerate}
        \item[(a)] $|M| = |M_1| \cup |M_2|,$

        \item[(b)] $ P^M = P^{M_1} \cup P^{M_2},$

        \item[(c)]  $ Q^M_\iota = Q^{M_1}_\iota \cup Q^{M_2}_\iota$ for $\iota = 0,1$.
    \end{enumerate}
    
    So $M$ is a $\tau$-model, moreover it is a model of $T$ as in $\oplus$ any 
    pair of distinct elements of $N^*_{\eta,\bar b}$  belongs to a relation, i.e. $\ell < k \le n \Rightarrow (a_\ell,a_k) \in P^M$.]
    
    \begin{enumerate}
        \item[$(*)_3$]  $T_{\eta,\bar k}$, the model completion of
        $T^0_{\eta,\bar k}$, is well defined and has elimination of quantifiers.
    \end{enumerate}
    
    [Why?  As $\tau$ is finite with no function symbols and $(*)_2$.]
    
    \begin{enumerate}
        \item[$(*)_4$]  $T_{\eta,\bar k}$ is $\NSOP_4$ (see \cite[2.5]{Sh:500}).
    \end{enumerate}
    
    [Why?  Because 
    
    \begin{enumerate}
    \item[$(*)_{4.1}$]   if (A) then (B), where:
    
    \begin{enumerate}
        \item[(A)]  
    
        \begin{enumerate}
            \item[(a)]  $ A_0,A_1,A_2,A_3$ are disjoint sets,
    
            \item[(b)] $ M_\ell$ is a model of  $T^0_{\eta,\bar k}$ with universe $A_\ell$ for $\ell = 0,1,2,3,$
    
            \item[(c)] if $\{\ell(1),\ell(2)\} \in \cW \coloneqq  \{\{0,1\},\{1,2\},\{2,3\},\{3,4\}\}$ then   $M_{\{\ell(1),\ell(2)\}}$ is a model of $T^0_{k,n}$ with universe $A_{\ell(1)} \cup A_{\ell(2)}$   extending $M_{\ell(1)}$ and $M_{\ell(2)}.$
        \end{enumerate}
        
        \item[(B)]  $M = \bigcup \{M_{\{\ell(1),\ell(2)\}} \colon \{\ell(1),\ell(2)\}
        \in \cW\}$ where the union is defined as in the proof of $(*)_2$, is a model of $T^0_{k,n}$ extending all of them.]
        \end{enumerate}
    \end{enumerate}
    
    [Why?  Clearly $M$ is a $\tau$-model and if $f$ embeds
    $N^*_{\eta,\bar k}$ into $M$, as in $(*)_2$ we have $\Rang(f) \subseteq
    M_{\ell(1), \ell(2)}$ for some $\{\ell(1),\ell(2)\} \in \cW$, a
    contradiction.]
    
    \begin{enumerate}
        \item[$(*)_5$]  $T_{\eta,\bar k}$ (and $\Mod_{T^0_{\eta,\bar k}}$) has the 
        $(\eta,\bar k)$-olive property as witnessed by $\varphi(x_0,x_1) = P(x_0,x_1),\psi_\iota(x_0,\dotsc,x_{k(\iota)}) = Q_\iota(x_0,\dotsc,x_{k(\iota)})$.
    \end{enumerate}
    
    [Why?  In Definition \ref{d2}(3), clause (a) holds trivially, and clause (c) is obvious from the choice
    of $T^0_{\eta,\bar k}$.  For clause (b)$_\lambda$ we are given $\langle f_\alpha \colon \alpha < \lambda\rangle$ where $f_\alpha$ is a function from $\alpha$ to $\{0,1\}$ and we have to find $M$ as there.  We define a $\tau$-model $M$ with:
    
    \begin{enumerate}
        \item[$\bullet$]  universe $\{a^*_\alpha \colon \alpha < \lambda\}$ such that $\alpha < \beta \Rightarrow a^*_\alpha \ne a^*_\beta,$
        
        \item[$\bullet$]  $P^M = \{(a^*_\alpha,a^*_\beta) \colon \alpha < \beta  < \lambda\},$  
        
        \item[$\bullet$]  $Q^M_\iota =\{(a^*_{\alpha_0},\dotsc,a^*_{\alpha_{k(\iota)-1}},a_\beta) \colon 
        \alpha_0 < \ldots < \alpha_{k(\iota)-1} < \beta$ 
        and $f_\beta \rest [\alpha_0,\alpha_{k(\iota)-1}]$ is constantly $\iota\}$.
    \end{enumerate}
    
    It suffices to prove that $M$ is a model of $T^0_{\eta,\bar k}$.  So 
    toward a contradiction assume $h$ embeds $N^*_{\eta,\bar k}$ 
    into $M$ and consider $h(a^*_\ell) =
    a_{g(\ell)}$ where $g \colon \{0,\dotsc,n\} \rightarrow \lambda$; recalling
    $\ell_1 < \ell_2 \le n \Rightarrow (a^*_{\ell_1} \ne a^*_{\ell_2})$
    and $h$ is an embedding, necessarily $g$ is a one-to-one function.  
    For $\ell < n$, recall that 
    $N^*_{\eta,\bar k} \models ``P(a^*_\ell,a^*_{\ell +1})"$ 
    but $h$ is an embedding so $M \models ``P[a^*_{g(\ell)},
    a^*_{g(\ell +1)}]"$, but if $g(\ell) \ge g(\ell +1)$ this fails by the choice of $P^M,$ hence $g(\ell) < g(\ell +1)$.  Now let $i_* \coloneqq
    \min\{i \colon \eta(i) = 1\}$.  Let $i_0 < \ldots <
    i_{k(0)-1}$ be from $\eta^{-1}\{0\}$ such that $i_0 = 0$ (recall
    Definition \ref{d2}(1)$(c)(\alpha)$ and $i_{k(\iota)-1}$ is maximal in
    $\eta^{-1}\{1\}$ hence $i_* \in [i_0,i_{k(0)-1})$.  Now
    $N^*_{\eta,\bar k} \models Q_0[a^*_{i_0},\dotsc,a^*_{i_{k(0)-1}},a^*_n]$  hence $M \models Q_0[a^*_{g(i_0)},\dotsc,
    a^*_{g(i_{k(0)-1})}, a^*_{g(n)}]$ and this implies that
    $f_{g(n)} \rest [g(i_0),g(i_{k(0)-1)})]$ is constantly $0$ hence in
    particular $f_{g(n)}(g(i_*))=0$.  Similarly let $j_0 < \ldots < j_{k(1)-1}$ be
    from $\eta^{-1}\{1\}$ such that $j_0 = i_*$; now $N^*_{\eta,\bar k} \models Q_1[a^*_{j_0},\dotsc, a^*_{j_{k(1)-1}},a^*_n]$ hence $M \models Q_1[a_{g(0)},\dotsc,a_{g(j_{k(1)-1})},a_n]$ hence, 
    $f_{g(n)} \, \rest \, 
    [g(j_0),g(j_{k(1)-1})]$ is constantly 1, hence in particular 
    $f_{g(n)}(g(i_*))=1$,
    a contradiction, so $(*)_5$ holds indeed.]
    
    \begin{enumerate}
        \item[$(*)_6$]  $T_{\eta,\bar k}$ has the $\SOP_3$.
    \end{enumerate}
    
    Why?  Again let $i_* = \min\{i \colon \eta(i)=1\},u_0 = \{0,\dotsc,i_*-1\},u_1
    = \{i_*\},$ and $u_2 = \{i_* +1,\dotsc,n\}$.
    
    Note that,
    
    \begin{enumerate}
        \item[$(*)_{6.1}$]  $(u_0,u_1,u_2)$ is a partition of $\{0,\dotsc,n\}.$
        
        \item[$(*)_{6.2}$]  if $\iota < 2$ and $(a_{\ell_0},\dotsc,a_{\ell_{k(\iota)}}) \in Q^{N^*_{\eta,\bar k}}_\iota$
        then $\{\ell_0,\dotsc,\ell_{k(\iota)}\} \cap u_j = \emptyset$ for some $j \le 2$.
    \end{enumerate}
    
    [Why?  Otherwise as $\ell_0 < \ell_1 < \ldots$ necessarily $\ell_0 <
    i_*,\ell_{k(\iota)} > \iota_*$ and $i_* = \ell_k$ for some $k \in
    (0,k(\iota))$.  But then $\eta(\ell_0) = 0 \ne 1 = \eta(k)$, in contradiction to Definition \ref{d2}(3)(c).]
    
    The rest should be clear by considering the proof of the model completion of the theory of triangle-free graphs having $\SOP_3$, see \cite[\S2]{Sh:500}.
\end{PROOF}

As in earlier cases, we apply a kind of guessing of clubs (almost suitable also for them i.e. for the proof with strict order and $\SOP_4$).   An unexpected gain is that here we use a weaker version:  there is no requirement $$\alpha < \lambda \Rightarrow \lambda > |\{C_\delta \cap \alpha \colon \delta \in S \text{ satisfies } \alpha \in \nacc(C_\delta)\}|,$$ but it is unclear how this helps.  Also here the use of the pair $(\bar{\cA},\bar{\mathbf g})$ may be helpful.

\begin{definition}\label{d13}
    1) For $\lambda$ regular uncountable and $\chi_2 > \chi_1 \ge \lambda$ let
    $\Qr_1(\chi_2,\chi_1,\lambda)$ mean that there are $S, \bar C,I,\bar{\cA},
    \bar{\mathbf g}$ witnessing it, which means that:
    
    \begin{enumerate}
        \item[$\boxplus$] 
    
        \begin{enumerate}
            \item[(a)] $S \subseteq \lambda$ and $I$ is an ideal on $S$,
    
            \item[(b)] $\bar C = \langle C_\delta \colon \delta \in S\rangle,$  
    
            \item[(c)]  $C_\delta \subseteq \delta$, note that possibly $\sup(C_\delta) < \delta,$
    
            \item[(d)] 
    
            \begin{enumerate}
                \item[$(\alpha)$]  $ \bar{\mathbf g} = \langle \bar g_j \colon j < \chi_2\rangle,$
    
                \item[$(\beta)$] $\bar g_j = \langle g_{j,\delta} \colon \delta \in S\rangle,$
    
                \item[$(\gamma)$] $g_{j,\delta} \colon C_\delta \rightarrow \{0,1\}.$ 
            \end{enumerate}
    
            \item[(e)] 
    
            \begin{enumerate}
                \item[$(\alpha)$]  $\bar{\cA} = \langle \bar{\cA}_j \colon j < \chi_2\rangle,$
    
                \item[$(\beta)$] $\bar{\cA}_j = \langle \cA_{j,\delta} \colon \delta \in S \rangle,$ 
    
                \item[$(\gamma)$] $\cA_{j,\delta} \subseteq \cP(\nacc(C_\delta).$
            \end{enumerate}
            
            \item[(f)] $\mathbf U_I(\chi_1) < \chi_2$; see Definition \ref{z20}, if $\chi_1=\lambda$ then we stipulate $\mathbf U_I(\chi_1) = \chi_1$ hence this means $\chi_1 < \chi_2,$
    
            \item[(g)] if $j_1 \ne j_2$ are $< \chi_2,\delta \in S, A_1 \in \cA_{j_1,\delta}$ and $A_2 \in \cA_{j_2,\delta},$ \then \, there is $\gamma \in A_1 \cap A_2$  such that $g_{j_1,\delta}(\gamma) \ne g_{j_2,\delta}(\gamma),$
    
            \item[(h)] if $j< \chi_2$ and $E$ is a club of $\lambda,$ \then \, for some $Y \in I^+$ hence $Y \subseteq S$  for every $\delta \in Y$ we have $\nacc(C_\delta) \cap E \in \cA_{j,\delta}$.
        \end{enumerate}
    \end{enumerate}
    
    2) For $\ell = 1,2,3$ let $\Qr_\ell(\chi_2, \chi_1, \lambda)$ be defined by:
    
    \begin{enumerate}
        \item[$\bullet$]  \If \, $\ell=1$ as above,
        
        \item[$\bullet$]  \If \, $\ell=2$ as above but there is a sequence $\langle J_\delta \colon \delta \in S\rangle$ of ideals on $\nacc(C_\delta)$ such that $\cA_{j, \delta} = \{\nacc(C_\delta) \backslash X \colon X \in J_\delta\},$ 
        
        \item[$\bullet$]  \If \, $\ell=3$ we use clauses (a)-(g) from part (1) and,

        \begin{enumerate}
            \item[(h)$^{-}$]  $(h)^- \quad$ if $E_j$ is a club of $\lambda$ for $j < \chi_2$ and $\langle \xi_j \colon j < \chi_2\rangle$ is a sequence of  ordinals with $\sup(\{\xi_j \colon j < \chi_2\}) < \chi_2$ \then \, we can find 
            $j_1 < j_2 < \chi_2, \delta \in S$ and $\gamma \in \nacc(C_\delta)$ such that   $\xi_{j_1} = \xi_{j_2},\gamma \in E_{j_1} \cap E_{j_2}$ and $g_{j_1,\delta}(\gamma) \ne g_{j_2,\delta}(\gamma)$.
        \end{enumerate} 
    \end{enumerate}
    
    3) $\Qr_{\ell,\iota}(\chi_2,\chi_1,\lambda)$ is defined as in $\Qr_\ell(\chi_2,\chi_1,\lambda)$ but $g_{j,\delta} \colon \nacc(C_\delta)
    \rightarrow \iota$, etc.
\end{definition}

\begin{remark}\label{d14}
    Can we weaken the conclusion of clause (h) of \ref{d13}(1), etc. to:
    
    \begin{enumerate}
        \item[$\bullet$]  $\{\alpha \in \nacc(C_\delta) \colon \sup(\alpha \cap E) > \max(C_\delta \cap \alpha)\} \in \cA_{j,\delta}$.
    \end{enumerate}
    
    That is, this suffices in Theorem \ref{d17} but there is no clear gain so we have not looked into it.
\end{remark}

\begin{fact}\label{d15}
    1) $\Qr_2(\chi_2,\chi_1,\lambda) \Rightarrow
    \Qr_1(\chi_2,\chi_1,\lambda) \Rightarrow \Qr_3(\chi_2,\chi_1,\lambda)$.
    
    2) We have $\Qr_1(\chi_2,\chi_1,\lambda)$ and even $\Qr_2(\chi_2,\chi_1,\lambda)$ \when:

    \begin{enumerate}
        \item[(a)] $\kappa^+ \leq \lambda \le \chi_1 < \chi_2 < 2^\kappa,$ 

        \item[(b)] $\kappa = \cf(\kappa),\lambda = \cf(\lambda),$ 

        \item[(c)] $\mathbf U_I(\chi_1) < \chi_2$ when $\lambda < \chi_1$ and $I$ an ideal on $S$ so $S \notin I,$

        \item[(d)]  $S \subseteq S^\lambda_\kappa$ is stationary, $\bar C = \langle C_\delta \colon  \delta \in S \rangle$ guess clubs, $C_\delta \subseteq \delta, \otp(C_\delta) = \kappa,$ 

        \item[(e)] $I = \{A \subseteq S$: for some club $E$ of $\lambda$  for no $\delta \in S$ do we have $\nacc(C_\delta) \cap E \in J^{\bd}_{C_\delta}\}$.
    \end{enumerate}
    
    3) If $\kappa^{+} < \lambda$ and clauses (a), (b) of part (2) hold and $S \subseteq S^\lambda_\kappa$ is stationary \then \, there is $\bar C$ as required in clause (d).
\end{fact}

\begin{PROOF}{\ref{d15}}  
    1) Easy.
    
    2) The proof is straightforward. 
    
    3)  Clause (d) follows by \cite[\S2]{Sh:420}. 
\end{PROOF}

\begin{theorem}\label{d17}
    1) If $T$ is complete, with the $3$-olive property  and $\lambda > \kappa^+$ and $\lambda,\kappa$ are regular,  $2^\kappa > \lambda \ge \kappa^{++} + |T|,$ \then \, $T$ has no universal model in $\lambda$ (for $\prec$).
    
    2) If $T$ is complete, with the $(\eta,\bar k,m)$-olive property  and $\lambda = \cf(\lambda) \ge |T|$ and $\Qr_1(\chi_2,\chi_1,\lambda),$ \then \, $\univ(\chi_1,\lambda,T) \ge \chi_2$.
    
    3) Similarly for class $\gk$ of models with the
    $\lambda-(\qf,\eta,\bar k)$-olive property,  see Definition \ref{d4}(2),  so e.g. for universal $K$ with the $\JEP$ and the  $\lambda-(\qf,\eta,\bar k)$-olive property.
    
    4) Like part (3) for a pair $\bar{\gk}$.
    
    5) We can weaken $\Qr_1(\chi_2,\chi_1,\lambda)$ to $\Qr_{1,\theta}(\chi_2,\chi_1,\lambda)$ \when \, $\theta = 2^\partial, \chi_1 = \chi^\partial_1,\partial < \lambda$.
\end{theorem}

\begin{remark}\label{d20}
    1) We can use $\Qr_3$ instead of $\Qr_1$ by the same proof but the gain is not clear.
    
    2) Assume $T$ is as in \ref{d17}(1), $\lambda \in \Univ(T)$. If e.g. $\lambda = \mu^+,\mu = \mu^{< \mu} = 2^\partial, \chi_1 = \lambda = \chi_2$ (so $T$ have a universal in $\lambda$), failure of $\Qr_1(\lambda,\lambda,\lambda)$ implies: there is $\cF \subseteq {}^\mu \mu$ such that $(\forall \eta \in {}^\mu \mu)(\exists \nu \in \cF)(\exists^\mu i < \mu)(\eta(i) = \nu(i))$.
\end{remark}

\begin{PROOF}{\ref{d17}}  
    1) It follows from (2) by Fact \ref{d15}(2).
    
    2) Let $(\bar\psi_0,\bar\psi_1)$, i.e. $\bar\psi_\iota = \langle  \psi_{\iota,k}(\bar x_0,\dotsc,\bar x_k) \colon k=1,\dotsc,k_\iota\rangle$
    for $\iota = 0,1$ witness the $(\Delta,\eta,\bar k,m)$-olive property. For simplicity we can, \wilog \, assume that $m=1$ and the formulas
    $\psi_{\iota}$ are quantifier-free and $T$ has only predicates  and its vocabulary is finite.  To make this proof also be a proof of Theorem
    \ref{d17}(3) let $\le_{\gk}$ be $\prec$ on $\mathrm{mod}_{T}$. Let $S,\bar C,\bar{\cA},\bar{\mathbf g}$ witness $\Qr_1(\chi_2,\chi_1,\lambda)$.  For each  $j < \chi_2$ we define $\bar f_j$ by:
    
    \begin{enumerate}
        \item[$(*)_1$]
    
        \begin{enumerate}
            \item[(a)] $\bar f_j = \langle f_{j,\alpha} \colon \alpha < \lambda \rangle,$
    
            \item[(b)]  $f_{j,\alpha} \colon \alpha \rightarrow \{0,1\}$ is defined by:
    
            \begin{enumerate}
                \item[$(\alpha)$] if $\beta < \alpha \in S$ then  $f_{j,\alpha}(\beta) = g_{j,\alpha}(\min(C_\alpha \setminus  \beta)),$
    
                \item[$(\beta)$]  if $\beta < \alpha \in \lambda \backslash S$ then $f_{j,\alpha}(\beta) = 0$.
            \end{enumerate}
        \end{enumerate}
    \end{enumerate}
    
    For each $j < \chi_2$ we can find $M_j
    \models T$ of cardinality $\lambda$ and pairwise distinct elements $\langle a_{j,\alpha} \colon \alpha < \lambda\rangle$ of $M_j$ satisfying
    (b)$_\lambda$ of Definition \ref{d2}(3) for $\bar f_j$.  Let $M_{j,\alpha} = M_j \rest \cup\{a_{j,\beta} \colon \beta < \alpha\}$.  
    Let the function  $h^0_j \colon \lambda \rightarrow M_j$ be defined by $h^0_j(\alpha) = a_{j,\alpha}$.
    
    Let $\cP \subseteq [\chi_1]^\lambda$ witness $\mathbf U_J(\chi_1)  < \chi_2$, so if $\lambda = \chi_1$ we use $\cP = \{\lambda\}$; \wilog
    \, $u \in \cP \wedge \alpha < \chi_1 \wedge |u \cap \alpha| = \lambda \Rightarrow u \cap \alpha \in \cP$.   For $u \in \cP$ or just $u \in [\chi_1]^\lambda$ let $h^1_u$ be a one-to-one function from $u$ onto $\lambda$.
    
    Towards a contradiction assume that there are 
    $\xi_* < \chi_2$ and a sequence
    $\langle \gA_\xi \colon \xi < \xi_*\rangle$ of models of $T$, each of
    cardinality $\le \chi_1$ witnessing $\univ(\chi_1,\lambda,T) < \chi_2$,
    even equal to $|\xi_*|$.  \Wilog \, the universe of $\gA_\xi$ is
    $\alpha_\xi \le \chi_1$ for $\xi < \xi_*$.  
    So for every $j < \chi_2$ there are $\xi =
    \xi_j < \xi_*$ and\footnote{Recall that now $\le_{\gk} \, = \, \prec \, \rest \, \mathrm{mod}_T$.}
    an $\le_{\gk}$-embedding $h^2_j$ of $M_j$ into
    $\gA_\xi$, hence there is $u_j \in \cP$ such that $W_j \coloneqq \{\alpha
    \in S \colon h^2_j(a_{j,\alpha}) \in u_j\} \in I^+$ and let $v_j \supseteq
    u_j \cup \Rang(h^2_j)$ be such that $v_j \in [\alpha_{\xi_j}]^\lambda
    \subseteq [\chi_1]^\lambda$ and 
    $\gA_j \rest v_j \prec \gA_j$ and let $\langle \gamma_{j,\alpha} \colon \alpha <
    \lambda\rangle$ list the members of $v_j$.
    
    Let $h'_j = h^2_{v_j} \circ (h^0_j \rest w_j)$ and let $h_j = h^1_{v_j} \circ h^2_j \circ (h^0_j \rest W_j)$; they are functions from $W_j$ into $\gA_{\xi_j}, \lambda$ respectively.  Let $N_j \coloneqq (\gA_{\xi_j} \rest v_j,P^{N_j}_*)$ be the expansion of $\gA_{\xi_j} \, \rest \, v_j$ by the relation $P^{N_j}_* = \Rang(h'_j)$ and
    let, 
    
    \begin{equation*}
        \begin{array}{clcr}
        E_j = \{\delta < \lambda \colon &\delta \text{ is a limit ordinal}, (\forall \alpha <
        \lambda)((h_{v_j})^{-1}(\alpha) \in  \{\gamma_{j,\beta} \colon \beta < \delta\} \\
        &\equiv \alpha < \delta) \text{ and } N_j \, \rest \,  \{\gamma_{j,\alpha} \colon \alpha <  \delta\} \prec N_j\},
        \end{array}
    \end{equation*}
    
    which clearly is a club in $\lambda$.   Hence by clause (h) of Definition \ref{d13}(1) there is an ordinal $\delta_j \in E_j \cap S$ such that $A_j  \coloneqq \nacc(C_\delta) \cap E_j$ belongs to $\cA_{j,\delta}$.
    
    As $\xi_* < \chi_2,|\cP| < \chi_2$ and $|\{h_j(a_{j,\delta}) \colon j < \chi_2,\delta \in S\}| < \sup\{\|\gA_\xi\| \colon \xi < \xi_*\} \le \chi_1  < \chi_2$ by the pigeon-hole-principle there are $j_1,j_2$ such that:
    
    \begin{enumerate}
        \item[$(*)_2$]
    
        \begin{enumerate}
            \item[(a)] $j_1 = j(1) < j_2 = j(2),$ 
    
            \item[(b)] $\xi_{j(1)} = \xi_{j(2)},$ 
    
            \item[(c)] $\delta_{j_1} = \delta_{j_2}$ call it $\delta$  (so $\delta \in S$), 
    
            \item[(d)] $u_{j_1} = u_{j_2}$ call it  $u$, so $u = u_{j_\iota} \subseteq |N_{j_\iota}|$ for $\iota = 1,2$,
    
            \item[(e)] $h^2_{j_1}(a_{j_1,\delta}) =  h^2_{j_2}(a_{j_2,\delta})$  call it $b$, so $b \in \Rang(h^2_{j_1}) \cap \Rang(h^2_{j_2})$. 
        \end{enumerate}
    \end{enumerate}
    
    By clause (g) of Definition \ref{d13}(1), there is  $\gamma \in A_{j_1} \cap A_{j_2}$ such that $g_{j_1,\delta}(\gamma) \ne g_{j_2,\delta}(\gamma)$.
    
    Now we shall choose $\alpha_\ell$ by induction on $\ell < n$ such that:
    
    \begin{enumerate}
        \item[$(*)_3$]
    
        \begin{enumerate}
            \item[(a)] $\alpha_\ell \in W_{j_{\eta(\ell)}},$ 
    
            \item[(b)] $\alpha_\ell < \gamma$ but $\alpha_\ell > \sup(C_\delta \cap \gamma),$ 
    
            \item[(c)] $\langle \alpha_0,\dotsc,\alpha_\ell \rangle$ is increasing, 
    
            \item[(d)]  in the model $N_{j_{\eta(\ell)+1}}$ the elements $h^2_{j_{\eta(\ell)+1}}(a_{j_1,\delta})=b, h^1_{j_{\eta(\ell)+1}}(a_{j_1,\alpha_\ell})$ realize the same quantifier type over $\{h_{j_{\eta(\ell(1))+1}}(a_{j_\eta,\alpha_{\ell(1)}})\colon \ell(1) < \ell\}$  or at least for all relevant (finitely many) formulas.
        \end{enumerate}
    \end{enumerate}
    
    If we succeed, then in the model $\gA_{\xi_*}$ which extends $N_{j_1}$ and $N_{j_2}$ the sequence $\langle h^2_{j_{\eta(\ell)}}(a_{j_{\eta(\ell)}, \alpha_\ell}) \colon \ell <
    n\rangle \char 94 \langle b \rangle$ realizes the ``forbidden" type, that is, the one from clause (c) of Definition \ref{d2}, which is a
    contradiction.
    
    As $\delta \in W_j \cap E_{j_{\eta(\ell)}}$ by the choice of $E_{j_{\eta(\ell)}}$ we can carry the induction.
    
    3) Similarly.
    
    4) As in \cite{Sh:457} and the above, just use $\partial$-tuples of $\bar a$'s.
\end{PROOF}

A sufficient condition for cases of $\Qr_i$ is the following.

\begin{definition}\label{d22}
    Let $\Qr_4(\lambda)$ mean: $\lambda = \mu^+$ and $\langle C_\delta,D_\delta \colon \delta \in S\rangle$ satisfies $C_\delta \subseteq
    \delta,D_\delta$ a filter on $\nacc(C_\delta)$ such that $\cP(\nacc(C_\delta)) \, / \, D_\delta$ satisfies the $2^\mu$-c.c. and, for every club $E$ of $\lambda$ for some $\delta \in S,E \cap \nacc(C_\delta) \in D^+_\delta$.
\end{definition}

Note the extreme case: 

\begin{conclusion}\label{z25}
    1) If (A) then (B), where: 

    \begin{enumerate}
        \item[(A)] 

        \begin{itemize}
            \item[$\bullet_{1}$] $\lambda = \mu^{+}$ and $2^{\mu} > \lambda \vee \gd_{\mu} = \lambda,$

            \item[$\bullet_{2}$] $\bar{C} = \langle C_{\delta} \colon \delta \in S \rangle$ guess club and $\otp(C_{\delta}) \geq \mu,$

            \item[$\bullet_{3}$] $T$ is a first-order complete theory with the olive property of cardinality $\leq \lambda.$
        \end{itemize}

        \item[(B)] $T$ has no universal model in $\lambda.$
    \end{enumerate}

    2) We can replace (A)$(\bullet_{1})$ by: 

    \begin{enumerate}
        \item[$\bullet$] $\gd_{\mu} > \lambda \wedge \cf(\mu) = \mu$ \underline{or just} there is $\cA \subseteq \{ C \subseteq \lambda \colon \otp(C) = \mu \}$ has cardinality $\lambda$ and it guess club of $\lambda.$  
    \end{enumerate}
\end{conclusion}

\begin{PROOF}{\ref{z25}}
    Easy.
\end{PROOF}

%%%%%%%%%%%%%%%%%%%%%%%%%%%%%%%%%%%%%%%
% Section 2
\newpage

\section{The class of Groups has the olive property} \label{2}

\subsection{General Groups}\label{A}\

We shall try to prove that the class of groups has a universal member almost only when cardinal arithmetic is close to G.C.H.  The following theorem does this.

\begin{theorem} \label{s2}
    The class of groups has the olive property (see Definition \ref{y36} or \ref{d2}(1)(d)$(\gamma)$), in fact, the  $(\eta,\bar k,m)$-olive property,  where $\eta = \langle 0,1,0 \rangle,  \bar k = (2,1),$ and $m=6$.
\end{theorem}

Why does Theorem \ref{s2} suffice?  Because then we can use Theorem  \ref{d17}(3); or see Conclusions \ref{s53}, \ref{s56}, we  break the proof into a series of definitions and claims; we may
replace the use of HNN extensions (in Claim \ref{s34}) and free amalgamation
(in Claim \ref{s31}) by the proof of Claim \ref{s50}.

\begin{definition}\label{s5}
    Let $\bar \psi \coloneqq \bar\psi^{\grp}_{\olive}$ be
    $(\psi_{0,1},\psi_{0,2}, \psi_{1,1})$ defined as follows (letting $m=6$):
    
    \begin{enumerate}
        \item[(a)]  $\psi_{0,1} = \varphi_0 = 
        \varphi_0(\bar x_{[m]},\bar y_{[m]}) = y^{-1}_5 x_0 y_5 = x_2,$
        
        \item[(b)]  $\psi_{1,1} = \varphi_1 = \varphi_1(\bar x_{[m]},
        \bar y_{[m]}) = x^{-1}_5 y_1 x_5 = y_3 \wedge x^{-1}_5 y_4 x_5 = y_4,$
        
        \item[(c)]  $\psi_{0,2} = \psi(\bar x_{[m]},\bar y_{[m]},\bar z_{[m]}) =
        (\sigma_*(x_0,y_1,z_4) = e \wedge \sigma_*(x_2,y_3,z_4) \ne e)$, on
        $\sigma_*,$ see below.
    \end{enumerate}
\end{definition}

\begin{dc}\label{s6}
    There is a $\sigma_* = \sigma_*(x,y,z)$ such that:
    
    \begin{enumerate}
    \item[(a)]  $\sigma_*$ is a group word,
    
    \item[(b)]  for some group $G$ and $a,b,c \in G$ we have ``$\sigma_*(a,b,c) \ne e_G$",
    
    \item[(c)]  for any group $G$ and $a,b,c \in G$ we have $e \in \{a,b,c\} \Rightarrow \sigma^G(a,b,c) = e_G.$
    \end{enumerate}
\end{dc}

\begin{PROOF}{\ref{s6}}
    Straightforward, e.g. $(x^{-1} y^{-1} x y)^{-1} 
    z^{-1}(x^{-1} y^{-1} xy)z$.
\end{PROOF}

\begin{claim}\label{s8}
    The $\bar\psi$ from Definition \ref{s5} satisfies clause (c) of Definition \ref{y36}, i.e., for no group $G$ and $\bar a_\ell \in {}^m G$ for $\ell < 4$ do the formulas
    there hold.
\end{claim}

\begin{remark}\label{s9}
    We prove more: there are no group $G$ and $\bar a_\ell \in {}^m G$ for
    $\ell=0,1,2,3$ such that $\varphi_0[\bar a_0,\bar a_1], \, \varphi_1
    [\bar a_1,\bar a_2], \, \varphi_1[\bar a_1,\bar a_3],$ and $\psi[\bar a_0,\bar
      a_2, \bar a_3]$.  
\end{remark}

\begin{PROOF}{\ref{s8}}
    Assume towards a contradiction that $G$ and $\langle \bar a_\ell \colon \ell < 4\rangle$ form a counterexample. Notice that conjugation by $a_{1,5}$ is an automorphism of $G,$ which we call $g$.
    
    Now, 
    
    \begin{enumerate}
        \item[$\bullet$]  $g(a_{0,0}) = a_{0,2}$ by  Definition \ref{s5}(a) as $G \models \varphi_0[\bar a_0,\bar a_1],$ 
        
        \item[$\bullet$]  $g(a_{2,1}) = a_{2,3}$ by first conjunct of Definition \ref{s5}(b) as $G \models \varphi_1[\bar a_1,\bar a_2],$
        
        \item[$\bullet$]  $g(a_{3,4}) = a_{3,4}$ by the second conjunct of Definition \ref{s5}(b) as $G \models \varphi_1[\bar a_1,\bar a_3]$.
    \end{enumerate}
    
    Together we have:
    
    \begin{enumerate}
        \item[$\bullet$]  $g(\sigma_*(a_{0,0},a_{2,1},a_{3,4})) =
         \sigma_*(a_{0,2},a_{2,3},a_{3,4}).$
    \end{enumerate}
    
    But this contradicts $G \models \psi[\bar a_0,\bar a_2,\bar a_3]$ (see clause Definition \ref{s5}(c)).
\end{PROOF}

\begin{definition}\label{s10}
    Let $\bar f \in \mathbf F_\lambda$, i.e. $\bar f = \langle f_\alpha \colon \alpha < \lambda\rangle $ where, for each $\alpha < \lambda, \,  f_\alpha \colon \alpha \rightarrow \{0,1\}$.
    
    1) Let $X_{\bar f}  \coloneqq  X_{\bar f,m},$ where we let $X_{\bar f,k}  \coloneqq \{x_{\alpha,\ell} \colon \alpha < \lambda,\ell < k\}$ for $k \le m$; recall
    that here $m=6$.
    
    2) Let $\bar x_{\alpha,k}  \coloneqq \langle x_{\alpha,\ell} \colon \ell < k\rangle$ for $k \le m$ and let $\bar x_\alpha  \coloneqq \bar x_{\alpha,m}$.

    3) For $\ell=0,1$ we define the set $\Gamma^\ell_{\bar f}$ of equations (pedantically, for $\ell=1$ conjunctions of two equations) as follows: 
    \[
    \{\varphi_\ell(\bar x_\alpha,\bar x_\beta) \colon \alpha < \beta < \lambda
    \text{ and } f_\beta(\alpha) = \ell\}.
    \]
    
    4) We define the set $\Gamma^2_{\bar f}$ of equations as follows: 
    \[
    \{\sigma_*(x_{\alpha,0},x_{\beta,1},x_{\gamma,4}) = e \colon \alpha < \beta <
    \gamma < \lambda \text{ and } f_\gamma \rest [\alpha,\beta] \text{ is
      constantly } 0\}.
    \]
    
    5) Let $G^5_{\bar f}$ be the group generated by $X_{\bar f,5}$ freely
    except for the equations in $\Gamma^2_{\bar f}$. Note that the $x_{\alpha,5}$'s are not mentioned in $\Gamma^2_{\bar f}$.
    
    6) Let $G^6_{\bar f}$ be the group generated by $X_{\bar f,6}$ freely except for the equations in $\Gamma^0_{\bar f} \cup \Gamma^1_{\bar f} \cup
    \Gamma^2_{\bar f}$.
\end{definition}

\begin{discussion}\label{s12}
    1) For our purpose we have to show that for $\alpha < \beta < \gamma$ (and $\bar f \in \mathbf F_\lambda$) we have: $$G_{\bar f,6} \models ``\psi[\bar x_\alpha,\bar x_\beta,\bar x_\gamma]" \text{ iff } f_\gamma \, \rest \,
    [\alpha,\beta] = 0_{[\alpha,\beta]}.$$ 
    
    For proving the ``if" implication, assume $f_\gamma \, \rest \, [\alpha,\beta] =
    0_{[\alpha,\beta]}$.  Now the satisfaction of
    ``$\sigma_*(x_{\alpha,0},x_{\beta,1},x_{\gamma,4}) = e"$
    is obvious by the role of $\Gamma^2_{\bar f}$, the analysis below is intended to prove the other half, ``$\sigma_*(x_{\alpha,2},x_{\beta,3},x_{\gamma,4}) \ne e$".  
    For proving the ``only if" implication it suffices to prove that
    ``$\sigma_*(x_{\alpha,0},x_{\beta,1},x_{\gamma,4}) \ne e"$ when $f_\gamma \rest [\alpha,\beta] \ne 0_{[\alpha,\beta]}$.  In both cases, we prove that this holds in $G^5_{\bar f}$ and then prove 
    that $G^5_{\bar f} \subseteq G^6_{\bar f}$ in a natural way.
    
    2) Of course, we also have to prove $G^6_{\bar f} \models \varphi_i(\bar x_\alpha,\bar x_\beta)$ when $\alpha < \beta < \lambda$
    and $f_\beta(\alpha) = \iota$.
\end{discussion}

\begin{claim}\label{s17}
    1) If $\alpha < \beta < \gamma < \lambda$ and $f_\gamma \rest [\alpha,\beta] \ne 0_{[\alpha,\beta]},$ \then \, $$G^5_{\bar f} \models
    ``\sigma_*[x_{\alpha,0},x_{\beta,1},x_{\gamma,4}] \ne e".$$ 
    
    2) If $\alpha < \beta < \gamma < \lambda,$ \then \,
    $G^5_{\bar f} \models ``\sigma_*(x_{\alpha,2},x_{\beta,3}, x_{\gamma,4}) \ne e"$. 
\end{claim}

\begin{PROOF}{\ref{s17}}
    1) Use \ref{s21} below with $X = \{x_{\xi,\ell} \colon \xi \in \{\alpha,\beta,\gamma\}$ and $\ell < 5\}$.
    
    2) Use \ref{s21}(2) below  with $X = \{x_{\xi,\ell} \colon \xi < \lambda,\ell < 5$ and $\ell > 0\}$.
\end{PROOF}

\begin{observation}\label{s21}
    1) If $x_{\alpha,\ell},x_{\beta,k} \in X^5_{\bar f}$ and $(\alpha,\ell) \ne (\beta,k),$ \then \, $G^5_{\bar f} \models ``x_{\alpha,\ell} \ne x_{\beta,k}"$.
    
    2) If $X \subseteq X^5_{\bar f}$ and
    $(\sigma_*(x_{\alpha,0},x_{\beta,1},
    x_{\gamma,4}) = e) \in \Gamma^2_{\bar f} \Rightarrow \{x_{\alpha,0},x_{\beta,1}, x_{\gamma,4}\} \nsubseteq X,$ \then \, $X$ generates freely a subgroup of $G^5_{\bar f}$.
\end{observation}

\begin{PROOF}{\ref{s21}}
    1) Let $G'  \coloneqq  \bigoplus \{\bbZ x \colon x \in X^5_{\bar f}\}$, it is an abelian
    group, and let $G''  \coloneqq \bigoplus \{\bbZ x_{\alpha, i} \colon \alpha < \lambda, i \notin
    \{\ell,k\}\}$, it is a subgroup. So $G'/G''$ by clause (c) of Definition
    \ref{s6} because $\{0,1,4\} \nsubseteq \{\ell, k\}$, satisfies all the equations in $\Gamma^2_{\bar f}$ and it satisfies the desired inequality.  As $G^5_{\bar f}$ is generated by
    $X^5_{\bar f}$ freely except for the equations in $\Gamma^2_{\bar f}$, the desired result follows.  Alternatively, use part (2).
    
    2) Let $H \coloneqq H_X$ be the group generated by $X$ freely.  We define a function $F$ from $X^5_{\bar f}$ into $H$ by: $$ F(x) \coloneqq 
    \begin{cases}
        x, \text{ if } x \in X, \\
        e_{H}, \text{ if } x \in X_{\bar{f}^{5}} \setminus X.
    \end{cases}$$
    
    Now $F$ respects every equation form $\Gamma^2_{\bar f}$ by  Claim \ref{s6}(c), hence $f$ induces a homomorphism from $G^5_{\bar f}$ into $H$, really onto.  Thus, the desired conclusion follows.
\end{PROOF}

\begin{definition}\label{s24}
    For $\beta < \lambda$ we define a partial function $F_\beta \colon X^5_{\bar f} \to X^5_{\bar f}$ as follows:
    
    \begin{enumerate}
        \item[$\bullet$]  if $\alpha < \beta$ and $f_\beta(\alpha) = 0,$ \then \, $F_\beta(x_{\alpha,0})  \coloneqq x_{\alpha,2},$
        
        \item[$\bullet$]  if $\gamma > \beta$ and $f_\gamma(\beta) = 1$ \then \, $F_\beta(x_{\gamma,1})  \coloneqq x_{\gamma,3}$ and $F_\beta(x_{\gamma,4})  \coloneqq x_{\gamma,4}$.
    \end{enumerate}
\end{definition}

\begin{claim}\label{s28}
    1) $F_\beta$ is a well-defined partial one-to-one 
    function from $X^5_{\bar f}$ to $X^5_{\bar f}$.
    
    2) The domain and the range of $F_\beta$ satisfy the criterion of Observation \ref{s21}(2). 
\end{claim}

\begin{PROOF}{\ref{s28}}
    1) It is a function as no $x_{\alpha,\ell}$ appears in two cases. Also if $F_\beta(x_{\alpha_1,\ell}) = x_{\alpha_2,k}$ then $\alpha_1 = \alpha_2 \wedge (\ell,k) \in \{(0,2),(1,3),(4,4)\},$ so $F_\beta$ is one-to-one.
    
    2) Assume $[\sigma_ (x_{\alpha_1,0},x_{\alpha_2,1},x_{\alpha_3,4})=e] \in \Gamma^2_{\bar f}$ so, 
    
    \begin{enumerate}
        \item[$(*)_1$]  $\alpha_1 < \alpha_2 < \alpha_3,$
    \end{enumerate}
    
    and, 
    
    \begin{enumerate}
        \item[$(*)_2$]   $f_{\alpha,3} \, \rest \,  [\alpha_1,\alpha_2] = 0_{[\alpha_1,\alpha_2]}$.
    \end{enumerate}
    
    First, toward contradiction assume
    $\{x_{\alpha_1,0},x_{\alpha_2,1},x_{\alpha_3,4}\} \subseteq \dom(F_\beta)$.  
    
    Now if $\alpha_1 \ge \beta$ then $x_{\alpha_1,0} \notin \dom(F_\beta)$, just inspect Definition \ref{s24} so necessarily $\alpha_1 < \beta$ and similarly $f_\beta(\alpha_1)=0$ (but not used).
    
    If $\alpha_2 \le \beta$ then $x_{\alpha_2,1} \notin \dom(F_\beta)$, so $\beta < \alpha_2$ and similarly $f_{\alpha_2}(\beta)=1$ (again not
    used) so together we get  $\alpha_1 < \beta < \alpha_2$.  Also as $x_{\alpha_3,4} \in \dom(F_\beta),$ it follows that ($\beta < \alpha_3$ which follows by earlier inequalities and) $f_{\alpha_3}(\beta)=1$, therefore $\beta$ witnesses that $f_{\alpha_3} \, \rest \,  [\alpha_1,\alpha_2]$ is not constantly zero; but
    this is a contradiction to $[\sigma(x_{\alpha_1,0},x_{\alpha_2,0},x_{\alpha_3,0}) = e] \in \Gamma^2_{\bar f}$.
    
    Second, assume towards  contradiction that
    $\{x_{\alpha_2,0},x_{\alpha_1,2},x_{\alpha_3,4}\} \subseteq  \Rang(F_\beta)$, but ``$x_{\alpha_2,0} \in \Rang(F_\beta)"$ is impossible by Definition \ref{s24}. 
\end{PROOF}

\begin{claim}\label{s31}
    To prove $G^5_{\bar f} \subseteq G^6_{\bar f}$
    any of the following conditions suffice:
    
    \begin{enumerate}
        \item[(a)]  there are a group $H$ extending $G^5_{\bar f}$ and $y_\zeta \in G$ for $\zeta < \lambda$ such that:  $$\zeta < \lambda \wedge F_\zeta(x_{\varepsilon_1,\ell_1}) = x_{\varepsilon_2,\ell_2} \Rightarrow H \models ``y^{-1}_\zeta x_{\varepsilon_1,\ell_1} y_\zeta
          = x_{\varepsilon_2,\ell_2}".$$
        
        \item[(b)]  for each $\zeta < \lambda$ there is a group $H$ extending $G^5_{\bar f}$ and $y \in G$ such that:  $$F_\zeta(x_{\varepsilon_1,\ell_1}) = x_{\varepsilon_2,\ell_2} \Rightarrow H \models ``y^{-1} x_{\varepsilon_1,\ell_1} y =
          x_{\varepsilon_2,\ell_2}".$$ 
    \end{enumerate}
\end{claim}

\begin{PROOF}{\ref{s31}}
    \underline{Clause (a) suffice}:

    We define a function $F$ from $X^6_{\bar f}$ into $H$ by: $$ F(x_{\varp, \ell}) \coloneqq 
    \begin{cases}
        x_{\varp, \ell}, \text{ if } \ell < 5 \wedge \varepsilon <
        \lambda, \\
        y_{\varp}, \text{ if }  \ell = 5  \wedge \varepsilon < \lambda,
    \end{cases}
    $$

    where in the first case,  $x_{\varp, \ell} \in G_{\bar{f}}^{5} \subseteq H.$ 
    
    Check that the mapping $F$ respects the equations in $\Gamma^0_{\bar f} \, \cup \, \Gamma^1_{\bar f} \, \cup \,  \Gamma^2_{\bar f}$ hence it induces a homomorphism $F^1$ from $G^6_{\bar f}$ into $H$, and  for every group word $$\sigma = \sigma(\ldots,x_{\varepsilon_i,\ell_i},\ldots)_{i<n}, x_{\varepsilon_i,\ell_i} \in X^5_{\bar f},$$ we have $G^6_{\bar f} \models ``\sigma = e" \Rightarrow G^5_{\bar f} \models ``\sigma = e"$, so we are done.
    
    \underline{Clause (b) suffice}:
    
    Let $(H_\zeta,y_\zeta)$ for $\zeta < \lambda$ 
    be as guaranteed by the assumption, i.e. clause (b).  \Wilog, \, $\zeta \ne \xi < \lambda = G_\zeta \cap G_\xi = G^5_{\bar f}$.  Now 
    clause (a) follows by using free amalgamation of
    $\langle H_\zeta \colon \zeta < \lambda\rangle$ over $G^5_{\bar f}$,  we know it is as required in clause (a), see e.g. \cite{LySc77}.
\end{PROOF}

\begin{claim}\label{s34}
    1) Clause (b) of Claim \ref{s31} holds.
    
    2) The conclusion of Claim \ref{s31} holds also for $G^6_{\bar f}$.
    
    3) The conclusions of Claim \ref{s17} hold also for $G^6_{\bar f}$.
\end{claim}

\begin{PROOF}{\ref{s34}}
    1) By the theorems on HNN extensions (see \cite{LySc77}) applied with the
    group being $G^5_{\bar f}$ and the partial automorphism $\pi_\zeta$ being the one $F_\zeta$ induced, i.e.,
    
    \begin{enumerate}
        \item[$\bullet$]  $\dom(\pi_\zeta)$ is the subgroup of $G^5_{\bar f}$ generated by $\dom(F_\zeta),$
        
        \item[$\bullet$]
        $\pi_\zeta(x_{\varepsilon,\ell}) =          F_\zeta(x_{\varepsilon,\ell})$ for $x_{\varepsilon,\ell} \in \dom(F_\zeta)$.
    \end{enumerate}
    
    By Claim \ref{s28}(2) and Observation \ref{s21}(2) we know that $\pi_\zeta$ is
    indeed an isomorphism.
    
    2) Follows by Claims \ref{s31} and \ref{s34}(1).
    
    3) By Claims \ref{s17} and \ref{s34}(2).
\end{PROOF}

Now, we prove Theorem \ref{s2}:

\begin{PROOF}{\ref{s2}}
    Should be clear by now.
\end{PROOF}

\subsection{Locally Finite Groups}\label{2B}

\begin{claim}\label{s50}
    The pair $(K_{\lfgr}, K_{\gr})$ of classes, i.e. 
    (locally finite groups, groups),
    has the olive property, as witnessed by $\bar\varphi$ from Definition \ref{s5}.
\end{claim}

\begin{PROOF}{\ref{s50}}
    We rely on Observation \ref{s51} below and use its notation.  Let $$J = \{(\alpha,\beta,\gamma) \colon \alpha < \beta < \gamma < \lambda \text{ and } f_\gamma \, \rest \,[\alpha,\beta]=0_{[\alpha,\beta]}\}.$$
    
    Let $G^5_{\bar f}, \, G^6_{\bar f}$ be as in the proof of Theorem \ref{s2}, that
    in Definition \ref{s10}.  Now
    for $\bar\alpha = (\alpha_0,\alpha_1,\alpha_2) \in J,$ let $\pi^5_{\bar\alpha}$ be the function from $X_{\bar f,5}$ (see Definition
    \ref{s10}(5)) into $K$, ($K$ is from \ref{s51}) defined as follows:
    
    \begin{enumerate}
        \item[$(*)_1$]  $\pi^5_{\bar\alpha}(x_{\beta,k})$ is
        
        \begin{enumerate}
            \item[$\bullet$]  $e_K,$ if $\beta \notin \{\alpha_0,\alpha_1,\alpha_2\},$
            
            \item[$\bullet$]  $z_{\ell,k},$ if $\beta = \alpha_\ell,\ell \le 2$.
        \end{enumerate}
    \end{enumerate}
    
    Now,
    
    \begin{enumerate}
        \item[$(*)_2$]  $\pi^5_{\bar\alpha}$ respects the equations from $\Gamma^2_{\bar f}$.
    \end{enumerate}
    
    [Why?  The equation $\sigma_*(x_{\alpha_0,0},x_{\alpha_1,1},x_{\alpha_2,4})=e$ holds as $K$ satisfies $\sigma_ (z_{0,0},z_{1,1},z_{2,4}) = e$.  For the other
    equations see Definition \ref{s6}(c); recall that the equations  for the cases of $\varphi_0,\varphi_1$ do not appear, see Definition  \ref{s10}(5).]
    
    Let $\pi^6_{\bar\alpha}$ be the following function from $X^6_{\bar f}$
    into $K$:
    
    \begin{enumerate}
        \item[$(*)_3$]  $\pi^6_{\bar\alpha}(x)$ is
        
        \begin{enumerate}
            \item[$\bullet$]  $\pi^5_{\bar\alpha}(x)$ when $x \in X^5_{\bar f},$
            
            \item[$\bullet$]  $z_s$ when $x=x_{\beta,5},\beta < \lambda$ and
            $s=s_{\bar\alpha, \beta} \coloneqq (\{\ell \le 2 \colon \alpha_\ell < \beta$ and $f_\beta(\alpha_\ell)=0\}, \{\ell \le 2 \colon \beta < \alpha_\ell$ and $f_{\alpha_\ell}(\beta)=1\})$.
        \end{enumerate}
    \end{enumerate}
    
    [Why is $\pi^6_{\bar\alpha}$ as required?  The least obvious point is: why $s \in S_*$?  Let $s = (u_1,u_2)$, now $\ell_1 \in u_1 \wedge \ell_2 \in u_2 \Rightarrow \alpha_{\ell_1} < \beta < \alpha_{\ell_2} \Rightarrow \ell_1 < \ell_2$ and $(\{0\},\{1,2\}) \ne s$ because $f_{\alpha_2} \, \rest \, [\alpha_0,\alpha_1]$ is constantly zero.]
    
    \begin{enumerate}
        \item[$(*)_4$]  $\pi^6_{\bar\alpha}$ respects the equations in $\Gamma^0_{\bar f} \cup \Gamma^1_{\bar f}$.
    \end{enumerate}
    
    [Why?  Check the definitions.]
    
    By $(*)_2,(*)_4$ there is a homomorphism $\pi_{\bar\alpha}$ from $G_{\bar f}$ into $K$ extending $\Pi^6_{\bar\alpha}$.  
    Let $G_*$ be the product of $J$-copies of $K$, i.e., 
    
    \begin{enumerate}
        \item[$(*)_5$]
        
        \begin{enumerate}
            \item[(a)] the set of elements of $G_*$ is the set of functions $g$ from $J$ into $K,$
    
            \item[(b)]  $G_* \models ``g_1 g_2 = g_3"$ \Iff \, $\bar\alpha \in J \Rightarrow K \models ``g_1(\bar\alpha) g_2(\bar\alpha) = g_3(\bar\alpha)"$.
        \end{enumerate}
    \end{enumerate}
    
    Now,
    
    \begin{enumerate}
        \item[$(*)_6$]  $G_*$ is a locally finite group.
        
        \item[$(*)_7$]  for $\alpha < \lambda, k<m$ let $\bar g_\beta = \langle g_{\beta,k} \colon k < m\rangle$ where $g^*_{\beta,k} \in G_*$ be defined by $(g_{\beta,k}(\bar\alpha))(x) = \pi^6_{\bar\alpha}(x_{\beta,k}).$
        
        \item[$(*)_8$]  $G_*,\langle \bar g_\beta \colon \beta < \lambda\rangle$ witnesses the olive property.
    \end{enumerate}
    
    [Why?  Check.]
    
    So we are done.
\end{PROOF}

\begin{observation}\label{s51}
    There are $K, z_{i,k}$ for $i < 3, \, k < m$ and $\langle \pi_s \colon s \in S_*\rangle$
    such that:
    
    \begin{enumerate}
        \item[(a)]  $K$ is a finite group,
        
        \item[(b)]  $z_{i,k} \in K,$
        
        \item[(c)]  $\sigma_*(z_{0,0},z_{1,2},z_{2,4})=e$ but $\sigma_*(z_{0,2},z_{1,3},z_{2,4}) \ne e,$ 
        
        \item[(d)]  $S_* = \{(u_1,u_2) \colon u_1,u_2 \subseteq \{0,1,2\}$ and $(\forall \ell_1 \in u_1)(\forall \ell_2 \in u_2)(\ell_1 < \ell_2)$ but $(u_1,u_2) \ne (\{0\},\{1,2\}),$
        
        \item[(e)]  for $s = (u_1,u_2) \in S_*$ we have: $\pi_s$ is a partial automorphism of $K$ such that:
        
        \begin{enumerate}
            \item[$(\alpha)$]  if $\ell \in u_1$ then $\pi_s(x_{\ell,0}) = x_{\ell,2},$
            
            \item[$(\beta)$]  if $\ell \in u_2$ then $\pi_s(x_{\ell,1}) = x_{\ell,2}, \pi_s(z_{\ell,4}) = z_{\ell,4},$
        \end{enumerate}
         
        \item[(f)]  moreover, there are $z_s \in K$ for $s \in S_*$ such that $(\forall x \in \dom(\pi_s))(\pi_s(x)=z^{-1}_s x z_s)$.
    \end{enumerate} 
\end{observation}

\begin{PROOF}{\ref{s51}}
    First, we ignore clause (f).  We use finite nilpotent groups.  Let $n_2  \coloneqq  6m, \, n_1  \coloneqq \binom{n_2}{2}, \, n_0 \coloneqq \binom{n_1}{2}$ and let $f_\ell \colon [n_{\ell +1}]^2 \rightarrow n_\ell$ be one-to-one for
    $\ell=0,1$.
    
    Let $K_1$ be the group generated by $\{y_{j,\ell} \colon j \le 2,\ell < 
    n_j\}$ freely except for the following equations:
    
    \begin{enumerate}
        \item[$(*)_1$]
        
        \begin{enumerate}
            \item[(a)] $y_{j,\ell} \cdot y_{j,\ell} = e,$
    
            \item[(b)] $ \quad [y_{j+1,\ell_1},y_{j+1,\ell_2}] = y_{j,f\{\ell_1,\ell_2\}}$, i.e. $y^{-1}_{j+1,\ell_1} y^{-1}_{j+1,\ell_2} y_{j+1,\ell_1} y_{j+1,\ell_2} = y_{j,f\{\ell_1,\ell_2\}},$ when $j<2,\ell_1 < \ell_2 < n_{j+1},$
     
            \item[(c)] $ [y_{j_1,\ell_1},y_{j_2,\ell_2}] = e$ when $(j_1 = 0 = j_2) \vee (j_1 \ne j_2 \le 2)$ and  $\ell_1 < n_{j_1},\ell_2 < n_{j_2}$.
        \end{enumerate}
    \end{enumerate}
    
    Clearly, $K_1$ is finite.
    
    Let $z'_{i,\ell} = y_{2,6i+\ell}$ for $i < 3,\ell < m$, let $\ell_*$ be
    such that $[[z'_{0,0},z'_{1,1}], \, z'_{2,4}] = y_{0,\ell_*}$.  Let $K_0$ be
    the subgroup $\{e,y_{0,\ell_*}\}$ of $K$, it is a normal subgroup as it is included in the center of $K_1$ and let $K_2  \coloneqq  K_1/K_0$ and we define $z_{i,\ell}$ as $z'_{i,\ell}/K_0$.
    
    Now,
    
    \begin{enumerate}
        \item[$(*)_2$]  $K_2,\langle z_{i,\ell} \colon i \le 2,\ell < m\rangle$ are as required in (a)-(e) of the claim.
    \end{enumerate}
    
    [Why?  We should just check that for $s \in S_*$ there is $\pi_s$ as required, i.e. that some subgroups of $K_2$ generated by subsets of $\langle z_{i,\ell} \colon i \le 2,\ell < m\rangle$ are isomorphic, but as none of them included $\{z_{0,0},z_{1,1},z_{2,4}\}$ and the way $K_2$ was defined this is straightforward.]
    
    Lastly, there is a finite group $K$ extending $K_2$ and $z_s \in K$ for $s \in S$ such that:
    
    \begin{enumerate}
        \item[$(*)_3$]   $x \in \dom(\pi_s) \Rightarrow z^{-1}_s x z_s = \pi_s(x)$.  
    \end{enumerate}
    
    Why?  Simply because $K_2$ can be considered as a group of permutations of the set $K_2$ (e.g. multiplying from the right), and it is easy to find $z_s \in \Sym(K_2)$ as required.
\end{PROOF}

\begin{conclusion}\label{s53}
    Assume $\Qr_1(\chi_1,\chi_2,\lambda)$.  
    
    Then there is no sequence $\langle G_\alpha \colon \alpha < \alpha_*\rangle$ of length $< \chi_2$ of groups of cardinality $\le \chi_1$ such that any locally finite group $H$ of cardinality $\lambda$ can be embedded into at least one of them.
\end{conclusion}

The following is an example.

\begin{conclusion}\label{s56}
    1) If $\mu = \cf(\mu),\mu^+ < \lambda = \cf(\lambda) < 2^\mu,$ \then \, there is no group of cardinality $\lambda$ universal for the class of locally finite groups.

    2) For example, if $\aleph_2 \le \lambda = 
    \cf(\lambda) < 2^{\aleph_0}$ this applies.
\end{conclusion}

\subsection{The Class of Groups is not Amenable}\label{2C}\

We have claimed (in earlier versions of \cite{Sh:820}) that the class of groups is amenable (see Dzamonja-Shelah \cite{Sh:614}) but this is not true.

An easy way to prove it is the following claim.

\begin{claim}\label{s60}
    For some group $G_0,$ if $G \supseteq G_0$ and $T = \Th(G)$ then $T$ is not
    amenable.
\end{claim}

\begin{PROOF}{\ref{s60}}
    In any forcing extension $\mathbf V_1 = \mathbf V^{\bbP}$ of $\mathbf V$  we
    have:
    
    \begin{enumerate}
        \item[$(*)$]   If $\mathbf V_1 \models ``\lambda  \coloneqq  \mu^+ = 2^\mu + \diamondsuit_{S^\lambda_\mu}, \, \bbQ$ is $\mu^+$-c.c.,  $(< \mu)$-complete forcing notion", \then \, in  $\mathbf V^{\bbQ}_1,T$ has no universal member in $\lambda$, moreover $\univ(\lambda,\lambda,T) \le \lambda^+$.
    \end{enumerate}
    
    [Why?  Because in $\mathbf V_1$ there is $\bar C = \langle C_\delta \colon \delta \in S^\lambda_\mu\rangle$ guessing clubs hence this
    holds also in $\mathbf V^{\bbQ}_1$, so we can apply Theorems \ref{d17}  and \ref{s2} or directly \ref{s53}.]  
    
    From $(*)$, by \cite{Sh:614} we get a contradiction to amenability (using a suitable $\bbQ$).
\end{PROOF}

%%%%%%%%%%%%%%%%%%%%%%%%%%%%%%%%%%%%%%%%%%
% References
\newpage

\bibliographystyle{amsalpha}
\bibliography{shlhetal}

\newcommand{\etalchar}[1]{$^{#1}$}
\providecommand{\bysame}{\leavevmode\hbox to3em{\hrulefill}\thinspace}
\providecommand{\MR}{\relax\ifhmode\unskip\space\fi MR }
\providecommand{\MRhref}[2]{%
  \href{http://www.ams.org/mathscinet-getitem?mr=#1}{#2}
}
\providecommand{\href}[2]{#2}
\begin{thebibliography}{She00b}

\bibitem[DS04]{Sh:614}
Mirna D{\v{z}}amonja and Saharon Shelah, \emph{{On the existence of universal
  models}}, Arch. Math. Logic \textbf{43} (2004), no.~7, 901--936,
  \href{https://arxiv.org/abs/math/9805149}{arXiv: math/9805149}. \MR{2096141}

\bibitem[D{\v{z}}a05]{Dj05}
Mirna D{\v{z}}amonja, \emph{Club guessing and the universal models}, Notre Dame
  J. Formal Logic \textbf{46} (2005), 283--300.

\bibitem[GS83]{Sh:174}
Rami~P. Grossberg and Saharon Shelah, \emph{{On universal locally finite
  groups}}, Israel J. Math. \textbf{44} (1983), no.~4, 289--302. \MR{710234}

\bibitem[KS92]{Sh:409}
Menachem Kojman and Saharon Shelah, \emph{{Nonexistence of universal orders in
  many cardinals}}, J. Symbolic Logic \textbf{57} (1992), no.~3, 875--891,
  \href{https://arxiv.org/abs/math/9209201}{arXiv: math/9209201}. \MR{1187454}

\bibitem[LS77]{LySc77}
Roger~C. Lyndon and Paul~E. Schupp, \emph{Combinatorial group theory},
  Ergebnisse der Mathematik und ihrer Grenzgebiete, vol.~89, Springer-Verlag,
  Berlin--Heidelberg--New York, 1977.

\bibitem[S{\etalchar{+}}a]{Sh:F1330}
S.~Shelah et~al., \emph{Tba}, In preparation. Preliminary number: Sh:F1330.

\bibitem[S{\etalchar{+}}b]{Sh:F1414}
\bysame, \emph{Tba}, In preparation. Preliminary number: Sh:F1414.

\bibitem[She80]{Sh:100}
Saharon Shelah, \emph{{Independence results}}, J. Symbolic Logic \textbf{45}
  (1980), no.~3, 563--573. \MR{583374}

\bibitem[She93a]{Sh:420}
\bysame, \emph{{Advances in cardinal arithmetic}}, {Finite and infinite
  combinatorics in sets and logic (Banff, AB, 1991)}, NATO Adv. Sci. Inst. Ser.
  C Math. Phys. Sci., vol. 411, Kluwer Acad. Publ., Dordrecht, 1993,
  \href{https://arxiv.org/abs/0708.1979}{arXiv: 0708.1979}, pp.~355--383.
  \MR{1261217}

\bibitem[She93b]{Sh:457}
\bysame, \emph{{The universality spectrum: consistency for more classes}},
  {Combinatorics, Paul Erd\H{o}s is eighty}, Bolyai Soc. Math. Stud., vol.~1,
  J\'anos Bolyai Math. Soc., Budapest, 1993,
  \href{https://arxiv.org/abs/math/9412229}{arXiv: math/9412229}, pp.~403--420.
  \MR{1249724}

\bibitem[She96]{Sh:500}
\bysame, \emph{{Toward classifying unstable theories}}, Ann. Pure Appl. Logic
  \textbf{80} (1996), no.~3, 229--255,
  \href{https://arxiv.org/abs/math/9508205}{arXiv: math/9508205}. \MR{1402297}

\bibitem[She00a]{Sh:589}
\bysame, \emph{{Applications of PCF theory}}, J. Symbolic Logic \textbf{65}
  (2000), no.~4, 1624--1674, \href{https://arxiv.org/abs/math/9804155}{arXiv:
  math/9804155}. \MR{1812172}

\bibitem[She00b]{Sh:702}
\bysame, \emph{{On what I do not understand (and have something to say), model
  theory}}, Math. Japon. \textbf{51} (2000), no.~2, 329--377,
  \href{https://arxiv.org/abs/math/9910158}{arXiv: math/9910158}. \MR{1747306}

\bibitem[She09]{Sh:h}
\bysame, \emph{{Classification theory for abstract elementary classes}},
  Studies in Logic (London), vol.~18, College Publications, London, 2009.
  \MR{2643267}

\bibitem[She17]{Sh:820}
\bysame, \emph{{Universal structures}}, Notre Dame J. Form. Log. \textbf{58}
  (2017), no.~2, 159--177, \href{https://arxiv.org/abs/math/0405159}{arXiv:
  math/0405159}. \MR{3634974}

\bibitem[SU06]{Sh:789}
Saharon Shelah and Alexander Usvyatsov, \emph{{Banach spaces and groups---order
  properties and universal models}}, Israel J. Math. \textbf{152} (2006),
  245--270, \href{https://arxiv.org/abs/math/0303325}{arXiv: math/0303325}.
  \MR{2214463}

\end{thebibliography}

\end{document}